\documentclass{article}

\usepackage{graphicx}              
\usepackage{amsmath}               
\usepackage{amsmath,scalerel}
\usepackage{amsfonts}              
\usepackage{amsthm}                
\usepackage{amssymb}
\usepackage{mathrsfs}
\usepackage[utf8]{inputenc}
\usepackage{calc}
\usepackage{hyperref}
\usepackage{multicol}
\usepackage{tikz-cd}
\usetikzlibrary{calc}
\usetikzlibrary{matrix}
\usepackage[margin=1in]{geometry} 
\usepackage{enumerate}
\usepackage{indentfirst}
\usepackage{float}
\usepackage[hypcap=false]{caption}
\newfloat{Diagram}{htbp}{loa}

\newtheorem{teo}{Theorem}
\newtheorem*{teo*}{Theorem}
\newtheorem{lem}[teo]{Lemma}
\newtheorem{prop}[teo]{Proposition}
\newtheorem{cor}[teo]{Corollary}
\newtheorem*{cor*}{Corollary}

\theoremstyle{definition}
\newtheorem*{ack}{Acknowledgements}

\newtheorem{obs}[teo]{Remark}
\newtheorem{defn}[teo]{Definition}

\makeatletter
\newcommand{\dobletilde}[1]{{%
 \text{\small$\mathpalette\double@widetilde{#1}$}%
}}
\newcommand{\double@widetilde}[2]{%
  \sbox\z@{$\m@th#1\widetilde{#2}$}%
  \ht\z@=.9\ht\z@
  \widetilde{\box\z@}%
}
\makeatother

\makeatletter
\newcommand{\minidobletilde}[1]{{%
 \text{\scriptsize$\mathpalette\double@widetilde{#1}$}%
}}

\makeatother

\newcommand{\dbtilde}[1]{{#1}''}
\newcommand{\cl}[1]{\mathbf{#1}} 
\newcommand{\CC}{\mathbb{C}}

\newcommand{\PP}{\mathbb{P}}
\newcommand{\ZZ}{\mathbb{Z}}      

\newcommand{\EE}{\mathbb{E}}

\newcommand\blfootnote[1]{%
  \begingroup
  \renewcommand\thefootnote{}\footnote{#1}%
  \addtocounter{footnote}{-1}%
  \endgroup
}
\newcommand{\eps}{\varepsilon}

\DeclareMathOperator{\LL}{\mathscr{L}}

\DeclareMathOperator{\im}{im}

\DeclareMathOperator{\rank}{rk}

\providecommand{\keywords}[1]
{
  \small	
  \textbf{\textit{Keywords---}} #1
}

\begin{document}

\nocite{*} 

\title{Degree of logarithmic foliations of type (1,1,1)}

\author{Mariano Chehebar$^1$}

\date{\vspace{-5ex}}

\maketitle

\begin{abstract} The space of codimension one holomorphic foliations of degree 1 in a projective space has an irreducible component whose general element is a logarithmic differential 1-form with simple poles in three hyperplanes. We compute its projective degree by resolving its rational parametrization map through succesive blow-ups with smooth centers.   
\end{abstract}
\blfootnote{\keywords{Foliation, Algebraic Geometry, Intersection Theory.}}

\section{Introduction}


A codimension one holomorphic foliation of degree $d$ in the complex projective space $\PP^n$ is defined by a non-zero twisted differential 1-form $\omega=\sum_{i=0}^n a_i dx_i$ such that $a_i$ are homogeneous polynomials of degree $d+1$ without non-constant common divisors, satisfying the descent condition $\sum_{i=0}^n a_ix_i=0$ and the Frobenius integrability condition $\omega\wedge d\omega=0$. Multiplication by a non-zero constant does not change the invariant hypersurfaces of $\omega$. Thus, the moduli space $\mathbb{F}(n,d)$ 
is defined as the Zariski closure of the set of codimension one holomorphic foliations of degree $d$ in $\PP^n$ inside the space $\PP H^0(\PP^n,\Omega_{\PP^n}^1(d+2))$ of projectivized global sections of twisted differential 1-forms. The geometry of the space of foliations in $\PP^n$ for $n\geq 3$ is very rich and inspired a lot of research. 

Many authors studied families and irreducible components of the space $\mathbb{F}(n,d)$. For general degree $d$, there are a few known irreducible components such as pullback components (see \cite{CELN}), rational components (see \cite{GMLN} and \cite{CPV}) and logarithmic components (studied in \cite{Omegar} and more recently in \cite{CGM}). A codimension one foliation belonging to a pullback component is defined by a 1-form which is the pullback by a rational map of a 1-form in $\PP^2$. To define logarithmic components, fix an integer partition $\cl{d}=(d_1,\ldots,d_m)$ of $d+2$ such that $m\geq 2$, $d+2=\sum_i d_i$ and $d_1\geq \cdots\geq d_m > 0$.  A logarithmic 1-form $\omega$ of type $\cl{d}$ in $\PP^n$ is defined by the formula $\omega=\left(\prod_{i=1}^m F_i\right)\sum_{i=1}^m\lambda_i\frac{dF_i}{F_i}$ where $F_i$ is a homogeneous polynomial of degree $d_i$ in the variables $x_0,\ldots,x_n$ and $\lambda_i\in\CC$ satisfy the descent condition $\sum_{i=1}^m \lambda_i d_i=0$. The Zariski closure of the space of logarithmic 1-forms of type $\cl{d}$ is called \textit{logarithmic component of type} $\cl{d}$. We denote it by $\mathscr{L}(\cl{d})$. When the partition is formed by two elements, i.e. $m=2$, we call it \textit{rational component} of type $(d_1,d_2)$. In \cite{Omegar}, O. Calvo Andrade proved for the first time that $\mathscr{L}(\cl{d})$ forms an irreducible component of the corresponding moduli space of foliations. Recently in \cite{CGM}, the authors gave an algebraic proof of this fact.

A related problem is to find the decomposition of $\mathbb{F}(n,d)$ into irreducible components for small $d$. In his book \cite{Jou}, Jouanolou exhausted these components for $d=0$ and $d=1$. In the case $d=0$, the space $\mathbb{F}(n,0)$ is irreducible. For $d=1$, the space $\mathbb{F}(n,1)$ (whose elements are also called \textit{Jacobi equations}) has two irreducible components: rational of type $(2,1)$ and logarithmic of type $(1,1,1)$. In their celebrated paper \cite{CLN}, Cerveau and Lins-Neto showed that there are six irreducible components of $\mathbb{F}(n,2)$, two of rational type, two of logarithmic type, one linear pullback component and another exceptional component, see \cite{CLN}. In the recent article \cite{CLP}, the authors studied the case $d=3$, showing that there are at least 24 irreducible components, but not exhausting them all.

An important geometric invariant of $\mathbb{F}(n,d)$ is the degree of its irreducible components. There have been several articles computing the degree of known irreducible components of the moduli space of foliations. In \cite{FV}, the authors gave formulas for the degree of some pullback components. The degree of the exceptional component has been computed in \cite{RV}. In \cite{CPV} and \cite{LV}, the authors found the degree of several rational components. The goal of this work is to
compute the degree of the logarithmic component of type $(1,1,1)$. For this, we resolve its rational parametrization map $\rho$ via three succesive blow-ups with smooth centers, see Diagram \ref{diag} at the end of Section 3:
\begin{teo*}
The map $\rho''':X'''\rightarrow \PP H^0(\PP^n,\Omega_{\PP^n}^1(3))$ depicted in Diagram \ref{diag} is a resolution of the rational parametrization map $\rho:X\dashrightarrow \PP H^0(\PP^n,\Omega_{\PP^n}^1(3))$ of the logarithmic component $\LL(1,1,1)$ of type $(1,1,1)$.
\end{teo*} 
We combine this Theorem with intersection theory techniques to get our main result:

\begin{teo*}
The degree of the logarithmic component $\LL(1,1,1)$ of type $(1,1,1)$ in $\PP^3$ is $80$.
\end{teo*}
This completes the degree computations for the irreducible components inside the moduli space $\mathbb{F}(n,1)$ of Jacobi equations in $\PP^n$:

\begin{cor*}
The moduli space $\mathbb{F}(3,1)$ of Jacobi equations in $\PP^3$ has two irreducible components $\LL(2,1)$ and $\LL(1,1,1)$ of degrees $55$ and $80$ respectively. 
\end{cor*}

In Section 2, we introduce logarithmic components and recall some useful results about them. Moreover, we describe the scheme structure of the base locus of the natural parametrization of the logarithmic component $\LL(1,1,1)$ and prove some techincal results for later use. In Section 3, we resolve the base locus of $\rho$ with a sequence of three succesive blow-ups with smooth centers. For this, we use local coordinates around a point of the exceptional divisor and compute the set-theoretical base locus in each of the three steps of the resolution. In Section 4, we use intersection theory techniques and the results of Sections 2 and 3 to compute the degree of the logarithmic component of type $(1,1,1)$. 

\begin{ack}
The author would like to thank Fernando Cukierman, Javier Gargiulo Acea, C\'esar Massri and Israel Vainsencher for their useful contributions, comments and suggestions.
\end{ack}

\section{Preliminaries}


In this section, we recall some basic definitions and facts about logarithmic components, see for example \cite{Omegar} and \cite{CGM}. Also, we study the scheme structure of the base locus of the natural parametrization of the logarithmic component of type $(1,1,1)$ and prove some technical results that we will use in the following sections.

\begin{defn}
Fix integer numbers $d\geq 0$, $n \geq 3$, $m\geq 2$ and an integer partition $\cl{d}=(d_1,\ldots,d_m)$ of $d+2$ such that $d+2=\sum_{i=1}^m d_i$ and $d_1\geq\cdots\geq d_m>0$. A logarithmic differential 1-form of type $\cl{d}$ in $\PP^n$ is defined by a twisted differential 1-form that can be written as $\omega=\left(\prod_{i=1}^m F_i\right)\sum_{i=1}^m\lambda_i\frac{dF_i}{F_i}$ where $F_i$ is a homogeneous polynomial of degree $d_i$ in the variables $x_0,\ldots,x_n$ and $\lambda_i\in\CC$ satisfy $\sum_{i=1}^m \lambda_i d_i=0$. The \textit{logarithmic component of type $\cl{d}$ in $\PP^n$} is the Zariski closure of the set of logarithmic differential 1-forms of type $\cl{d}$ in $\PP^n$. We denote it by $\LL(\cl{d})$.
\end{defn}

In other words, there is a multilinear map
\begin{align*}
\mu_{\cl{d}}:\Lambda(\cl{d})\times \prod_{i=1}^m S_n(d_i)\rightarrow H^0(\PP^n,\Omega_{\PP^n}^1(d+2))
\end{align*}
 defined by $\mu_{\cl{d}}(\cl{\lambda},\cl{F})= F\sum_{i=1}^m\lambda_i\frac{dF_i}{F_i}$ where $S_n(d_i)=H^0(\PP^n,\mathcal{O}(d_i))$, $\cl{F}=(F_1,\ldots,F_m)\in \prod_{i=1}^m S_n(d_i)$,  $F=\prod_{i=1}^m F_i$ and $\Lambda(\cl{d})=\left\{\cl{\lambda}=(\lambda_1,\ldots,\lambda_m)\in\CC^m:\sum_{i=1}^m \lambda_i d_i=0\right\}$. This multilinear map induces a rational map
\begin{align*}
\rho_{\cl{d}}:X_{\cl{d}}:=\PP(\Lambda(\cl{d}))\times \prod_{i=1}^m \PP S_n(d_i)\dashrightarrow \PP H^0(\PP^n,\Omega_{\PP^n}^1(d+2)).
\end{align*}
The logarithmic component $\mathscr{L}(\cl{d})$ of type $\cl{d}$ is the Zariski closure of the image of the map  $\rho_{\cl{d}}$. Notice that the target space of $\rho_{\cl{d}}$ is a projective space. Thus, $\rho_{\cl{d}}$ is defined by a list of global sections of the line bundle $\mathcal{O}_{X_{\cl{d}}}(1,\ldots,1)$. To avoid overloading notation, we will usually omit projective classes for elements in the image of $\rho_{\cl{d}}$, writing for example $\rho_{\cl{d}}([(\cl{\lambda},\cl{F})])= F\sum_{i=1}^m\lambda_i\frac{dF_i}{F_i}$.

\begin{obs}
Notice that every logarithmic differential 1-form $\omega=F\sum_{i=1}^m\lambda_i\frac{dF_i}{F_i}$ of type $\cl{d}$ in $\PP^n$ defines a codimension one foliation of degree $d$. Indeed, straightforward calculations yield the two equalities $i_R(\omega)=\left(\sum_{i=1}^m \lambda_i d_i\right)F$ and $d\omega=\frac{dF}{F}\wedge \omega$. Therefore $i_R(\omega)=0$ and $\omega\wedge d\omega=0$.
\end{obs}

As its name already suggests, the varieties $\mathscr{L}(\cl{d})$ are irreducible components of the space of codimension one foliations.

\begin{teo}[\cite{Omegar} and \cite{CGM}] The variety $\mathscr{L}(\cl{d})\subseteq \mathbb{F}(n,d)$ of logarithmic forms of type $\cl{d}$ is an irreducible component of $\mathbb{F}(n,d)$. Moreover, the scheme $\mathbb{F}(n,d)$ is generically reduced along $\mathscr{L}(\cl{d})$. 
\end{teo}

Even though the map $\rho_{\cl{d}}$ is not injective in general, it is generically finite whenever $m\geq 3$.  

\begin{obs} \label{defS(d)}
Notice that $\rho_{\cl{d}}$ cannot be generically injective in general. Indeed, fix an integer partition $\cl{d}$ as before and define $A_e=\left\{i:d_i=e\right\}$. Let $\mathbb{S}(e)=\left\{\sigma\in\mathbb{S}_m:\sigma(j)=j\text{ for all } j\notin A_e\right\}$ be a subgroup of the group of permutations of $m$ elements and let $\mathbb{S}(\cl{d})=\prod_e\mathbb{S}(e)\subseteq\mathbb{S}_m$. Notice that $\mathbb{S}(\cl{d})$ acts in the domain of $\mu_{\cl{d}}$ via $\sigma\cdot (\cl{\lambda},\cl{F})=(\sigma\cdot\cl{\lambda},\sigma\cdot\cl{F})$ in the natural way. Thus, the map $\mu_{\cl{d}}$ is constant in the orbits of this action.    
\end{obs} 

\noindent
We say that $\rho_{\cl{d}}$ is injective \textit{up to order} if the induced map with domain $\left(\Lambda(\cl{d})\times \prod_{i=1}^m S_n(d_i)\right)/\mathbb{S}(\cl{d})$ is injective.
 
 \begin{prop}[\cite{CGM}] \label{geninyupto}
 Let $\cl{d}=(d_1,\ldots,d_m)$ be an integer partition with $d_1\geq\cdots\geq d_m>0$ and $m\geq 3$. The rational map $\rho_{\cl{d}}$ is generically injective up to order. 
 \end{prop}

For computing the degree of a logarithmic component, it is essential to understand its base locus. The set theoretical base locus of the natural parametrization $\rho_{\cl{d}}$ has been described in \cite[Section 5]{CGM}. Let us recall briefly that description. Let $\cl{d}$ be as above and take a decomposition $\varphi=(e,\cl{d'})\in\ZZ_{\geq 0}^{m\times m'}\times (\ZZ_{\geq 0}^{m'}\setminus\vec{0})$ such that $d_i=\sum_{j=1}^{m'}e_{ij}d'_j$ for every $1\leq i\leq m$. For each one of those decompositions, define the Segre-Veronese map 
\begin{align*}
\nu_{\varphi}:\prod_{j=1}^{m'} S_n(d'_j)\dashrightarrow \prod_{j=1}^{m} S_n(d_j)
\end{align*}
given by $\nu_{\varphi}(G_1,\ldots,G_{m'})=(\prod_{j=1}^{m'}G_j^{e_{1j}},\ldots,\prod_{j=1}^{m'}G_j^{e_{mj}})$. Also, define sets $B_{\varphi}:=\Lambda(\varphi)\times \im(\nu_{\varphi})$ where $\Lambda(\varphi)=\left\{\cl{\lambda}=(\lambda_1,\ldots,\lambda_m)\in\CC^m:\sum_{i=1}^m \lambda_i e_{ij}=0 \text{ for all }j\right\}$. Notice that the sets $B_{\varphi}$ are contained in the base locus of $\mu_{\cl{d}}$ and therefore their projectivizations are included in the base locus of $\rho_{\cl{d}}$. Also, for every pair of decompositions $\varphi^{(1)}=(e^{(1)},\cl{d'^{(1)}}),\varphi^{(2)}=(e^{(2)},\cl{d'^{(2)}})$, define a relation in the following way: $\varphi^{(1)}\leq\varphi^{(2)}$ if $\text{rank}(e^{(1)})=\text{rank}(e^{(2)})$ and there exists a matrix $e^{(3)}$ such that $e^{(1)}=e^{(2)}\cdot e^{(3)}$. This realtion is not a partial order. Indeed, if two matrices $e^{(1)},e^{(2)}$ differ by a permutation of columns, the corresponding decompositions $\varphi^{(1)},\varphi^{(2)}$ satisfy $\varphi^{(1)}\leq \varphi^{(2)}$ and $\varphi^{(1)}\geq \varphi^{(2)}$ at the same time. Therefore, we will consider decompositions up to permutation of columns so that $\leq$ is an order relation.      
The irreducible components of this base locus are in one-to-one correspondence with vector partitions that are maximal with respect to that order.

\begin{teo}[\cite{CGM}]
The irreducible components of the base locus of $\rho_{\cl{d}}$ are the projectivization of the sets $B_{\varphi}$ where $\varphi=(e,\cl{d'})$ is a maximal decomposition with respect to the order $\leq$ and $\rank(e)< m$. 
\end{teo}

As we are going to focus in the case $d=1,m=3$ and $\cl{d}=(1,1,1)$, let us set the notation $\rho:=\rho_{(1,1,1)}$, $\Lambda:=\Lambda(1,1,1)$ and $X:=X_{(1,1,1)}\simeq \PP^1\times (\PP^n)^3$. Also, let $B_{\rho_{\cl{d}}}$ be the base locus scheme of $\rho_{\cl{d}}$, let $B_{\rho}$ be the base locus scheme of $\rho$ and set $\cl{\mu^{(1)}}=(0,1,-1), \cl{\mu^{(2)}}=(1,0,-1), \cl{\mu^{(3)}}=(1,-1,0)\in \Lambda $. A direct application of the former result to our special case yields the decomposition of the set-theoretical base locus of $\rho$ into the following four irreducible components.

\begin{cor}
The irreducible components of the base locus of $\rho$ are given set theoretically by the following smooth varieties:
\begin{align*}
(B_0)_{\text{red}}=&\left\{[(\cl{\lambda},\cl{F})]\in X:F_1=F_2=F_3 \right\} \\
B_1=&\left\{[(\cl{\mu^{(1)}},\cl{F})]\in X:F_2=F_3 \right\} \\
B_2=&\left\{[(\cl{\mu^{(2)}},\cl{F})]\in X:F_1=F_3 \right\} \\
B_3=&\left\{[(\cl{\mu^{(3)}},\cl{F})]\in X:F_1=F_2 \right\}.
\end{align*}
\end{cor}

\begin{obs}
The components $B_1,B_2$ and $B_3$ are smooth of dimension $2n$ and $(B_0)_{\text{red}}$ is also smooth and its dimension is $n+1$. Additionally, the sets $B_1$, $B_2$ and $B_3$ are pairwise disjoint and the intersections $(B_0)_{\text{red}}\cap B_i=\left\{[(\cl{\mu^{(i)}},\cl{F})]\in X:F_1=F_2=F_3 \right\}$ are smooth of dimension $n$ for $i=1,2,3$.
\end{obs}

We will show that the base locus scheme $B_{\rho}$ of $\rho$ is reduced along $B_1,B_2,B_3$ and non-reduced along $(B_0)_{\text{red}}$. For this, let us describe the Zariski tangent space of the base locus scheme $B_{\rho}$ along its irreducible components. Let $W$ be a vector space and for every $w\in W$ let $[w]$ be its class in $\PP W$. There is a natural identification of the tangent space $T_{[w]} \PP W$ of $\PP W$ at $[w]$ with $W/[w]$. From now on, let $(\cl{\lambda'},\cl{F'})$ be an element of the vector space $\Lambda\times(S_n(1))^3$ and let $[(\cl{\lambda'},\cl{F'})]_x$ be its class inside $T_xX$ for some $x\in X$. We will omit the subindex $x$ whenever it is clear from the context. In order to compute the tangent space of $B_{\rho}$, we prove the following techincal result.

\begin{lem} \label{lema1}
Let $\cl{d}=(d_1,\ldots,d_m)$ be a partition of an integer, let $x=(\cl{\lambda},\cl{F}),v=(\cl{\lambda'},\cl{F'})$ be elements of $\Lambda(\cl{d})\times \prod_{i=1}^m S_n(d_i)$ such that the class of $x$ belongs to the base locus $B_{\rho_{\cl{d}}}$ of $\rho_{\cl{d}}$ and let $\eps$ be a formal parameter. 
The following formula holds:
\begin{align} \label{flalema1}
\rho([x+\eps v])= \eps H_1(x,v)+\eps^2 H_2(x,v) \text{ mod}(\eps^3),
\end{align}
where the terms in the formula are defined by 
\begin{align*}
H_1(x,v)=F\left(\sum_{i=1}^m \lambda'_i \frac{dF_i}{F_i}+d\left(\sum_{i=1}^m\lambda_i\frac{ F'_i}{F_i}\right)\right) \\
H_2(x,v)=\left(\sum_{i=1}^m \frac{F'_i}{F_i}\right)H_1(x,v)+F d\left(\sum_{i=1}^m\frac{ \lambda'_i F'_i F_i-\frac{\lambda_i}{2}{F'_i}^2}{{F_i}^2}\right).
\end{align*} 
 
\begin{proof}
Notice that $\frac{1}{F_i+\eps F'_i}=\frac{1}{F_i(1+\eps(\frac{F'_i}{F_i}))}\equiv \frac{1}{F_i}\left(1-\eps\frac{F'_i}{F_i}+\eps^2\left(\frac{F'_i}{F_i}\right)^2\right) \text{ mod}(\eps^3)$. Therefore, we replace in the formula of $\rho$ to get
\begin{align*}
\rho([x+\eps v])\equiv \left(\prod_{i=1}^m (F_i+\eps F'_i)\right)\left(\sum_{i=1}^m (\lambda_i+\eps\lambda'_i)\frac{d(F_i+\eps F'_i)}{F_i}\left(1-\eps\frac{F'_i}{F_i}+\eps^2\left(\frac{F'_i}{F_i}\right)^2\right)\right) \text{ mod}(\eps^3).
\end{align*}  
From a direct expansion of the product on the left and using the hypothesis $\sum_{i=1}^m \lambda_i\frac{dF_i}{F_i}=0$, we get 
\begin{align*}
\rho([x+\eps v])\equiv F\left(\sum_{i=1}^m \lambda'_i \frac{dF_i}{F_i}+\sum_{i=1}^m \lambda_i\left(\frac{ dF'_i}{F_i}-\frac{ F'_idF_i}{F_i^2}\right)\right)\left(\eps+ \left(\sum_{i=1}^m \frac{F'_i}{F_i}\right)\eps^2\right) + \\
+F \left(\sum_{i=1}^m \lambda'_i \left(\frac{ dF'_i}{F_i}-\frac{ F'_idF_i}{F_i^2}\right)-\lambda_i\left(\frac{F'_i dF'_i}{F_i^2}-\frac{dF_i{F'_i}^2}{F_i^3}\right)\right) \eps^2 \text{ mod}(\eps^3).
\end{align*}
We replace the equalities $\frac{ dF'_i}{F_i}-\frac{ F'_idF_i}{F_i^2}=d\left(\frac{F'_i}{F_i}\right)$ and $\frac{F'_i dF'_i}{F_i^2}-\frac{dF_i{F'_i}^2}{F_i^3}=\frac{1}{2}d\left(\frac{{F'_i}^2}{F_i^2}\right)$ in the last formula. The result follows.
\end{proof} 
 
\end{lem}

\begin{obs} 
Notice that the expression $H_1(x,v)$ in \eqref{flalema1} does not depend on the representative of a class in $T_xX$, as the point $x$ is in the base locus of $\rho$. Therefore, the expression $H_1$ coincides with the derivative of the parametrization $d\mu_{\cl{d}}$ in the closed points of the base locus, see \cite[Section 6]{CGM}.
\end{obs}

We use the previous Lemma to get information about the scheme structure of the base locus of $\rho$.

\begin{prop} \label{propbl1}
The base locus scheme $B_{\rho}$ of $\rho$ is reduced along $B_1$, $B_2$ and $B_3$ and non-reduced along $(B_0)_{\text{red}}$.

\begin{proof}
Let us compute the Zariski tangent space of $B_{\rho}$ at a general closed point of each irreducible component. Take a closed point $[(\cl{\lambda},\cl{F})]\in B_{\rho}$ and a tangent vector $[(\cl{\lambda'},\cl{F'})]\in T_x X\simeq  \left(\Lambda/\langle\cl{\lambda}\right)\rangle\times \prod_{i=1}^3 \left(S_n(1)/\langle F_i\rangle\right) $. The vector $[(\cl{\lambda'},\cl{F'})]$ is tangent to $B_{\rho}$ at $[(\cl{\lambda},\cl{F})]$ if and only if $\rho((\cl{\lambda},\cl{F})+\eps(\cl{\lambda'},\cl{F'}))=0 \text{ mod} (\eps^2)$.
Using Lemma \ref{lema1} and dividing by $F$, we get the condition
\begin{align} \label{eq2}
\frac{H_1((\cl{\lambda},\cl{F}),(\cl{\lambda'},\cl{F'}))}{F}=\sum_{i=1}^3 \lambda'_i\frac{dF_i}{F_i}+d\left(\sum_{i=1}^3 \lambda_i\frac{F'_i}{F_i}\right)=0.
\end{align}

If we assume $F_1=F_2=F_3$ holds, equality \eqref{eq2} becomes $d\left(\sum_{i=1}^3 \lambda_i\frac{F'_i}{F_1}\right)=0$. Thus, the tangent space of $B_{\rho}$ along a point $[(\cl{\lambda},\cl{F})]\in (B_0)_{\text{red}}$ is defined by
\begin{align} \label{eq3}
\sum_{i=1}^3 \lambda_i F'_i=0 \text{ mod}(F_1)
\end{align}
and therefore has codimension $n$ inside $T_{[(\cl{\lambda},\cl{F})]}X$. We deduce that $B_{\rho}$ is non-reduced along $(B_0)_{\text{red}}$.

For the case $F_1=F_2\neq F_3$ and $\cl{\lambda}=(1,-1,0)$, equality \eqref{eq2} becomes $(\lambda'_1+\lambda'_2)\frac{dF_1}{F_1}+\lambda_3'\frac{dF_3}{F_3}+d\left(\frac{F'_1-F'_2}{F_1}\right)=0$. As $\lambda'_1+\lambda'_2=-\lambda'_3$, the tangent space of $B_{\rho}$ along a point $(\cl{\lambda},\cl{F})\in (B_0)_{\text{red}}$ is defined by the two conditions
\begin{align*} 
 F'_1=F'_2 \text{ mod} (F_1) & & \lambda'_3=0. 
\end{align*}
Therefore $B_{\rho}$ is reduced along $B_3$ as its tangent space has codimension $n+1$ inside $X$. The argument for reducedness along $B_1$ and $B_2$ is analogous. 
\end{proof}

\end{prop}

The reader should note that we computed the Zariski tangent space of $B_{\rho}$ along the non-reduced component supported in $(B_0)_{\text{red}}$ in the proof of the last proposition, see equation \eqref{eq3}. 

\begin{defn}
We define $B_0$ to be the scheme-theoretical (non-reduced) irreducible component of $B_{\rho}$ supported in $(B_0)_{\text{red}}$.
\end{defn}

To finish the section, we point out a final remark regarding the tangent space of the scheme $B_0$.

\begin{obs} \label{secB0}
 From the proof of Proposition \ref{propbl1}, we obtain that $TB_{0} |_{(B_0)_{\text{red}}}$ is a vector bundle of rank $2n+1$ and $T_x B_0=T_x B_i+T_x (B_0)_{\text{red}}$ for every $x\in B_i\cap (B_0)_{\text{red}}$, $i=1,2,3$. Also, using the obvious isomorphism $(B_0)_{\text{red}}\simeq \PP\Lambda\times\PP S_n(1)$ and the description of the Zariski tangent space of $B_{\rho}$, we deduce the following exact sequence of vector bundles on $(B_0)_{\text{red}}$:
\[0 \rightarrow TB_{0} |_{(B_0)_{\text{red}}} \rightarrow TX|_{(B_0)_{\text{red}}} \rightarrow T\PP S_n(1) \otimes \mathscr{O}_{\PP \Lambda}(1) \rightarrow 0.\]	
Here the surjective map from the right is given on each fiber over $(\lambda,(F_0,F_0,F_0))\in (B_0)_{\text{red}}$ by the assignment $[(\cl{\lambda'},\cl{F'})] \mapsto [\sum_i \lambda_i F'_i]$. We observe that $[(\lambda',0)]$ represents a tangent vector at $[(\cl{\lambda},(F_0,F_0,F_0))]$ to $(B_0)_{\text{red}}$ and $B_{\rho}$ as well.
\end{obs}


\section{Resolution of the rational parametrization map}

In this section, we resolve the indeterminacy locus of $\rho$ through a sequence of three succesive blow-ups with smooth centers. 
In order to do this, we use standard local coordinates of a blow-up in an open neighbourhood of a point of the exceptional divisor. We compute the support of the base locus in each of the three steps of the resolution. There is a commutative diagram summarizing the three steps of the resolution of $\rho$ at the end of the section, see Diagram \ref{diag}. From now on, the reader is advised to use it as an outline for organizing the whole argument.

Let us start with the first step of the resolution of $\rho$.

\begin{defn}
Denote the blow-up of $X$ along $(B_0)_{\text{red}}$ by $X'$, its natural projection map by $\pi:X'\rightarrow X$, its exceptional divisor by $E'\simeq \PP N_{(B_0)_{\text{red}}} X$, the strict transforms of $B_i$ for $i=1,2,3$ by $B_i'$ and the extension of the rational map $\rho$ by $\rho':X'\dashrightarrow \PP H^0(\PP^n,\Omega_{\PP^n}^1(3))$. The map $\rho'$ is well defined on an open subset of $E'$ and coincides  with the rational map $\rho \circ \pi$ on $X'-E'$. Let $B_{\rho'}$ be its indeterminacy scheme and let $B'_0$ be the closed subscheme of $B_{\rho'}$ formed by all the scheme-theoretical (not necessarily reduced) irreducible components of $B_{\rho'}$ supported in $E'$. 
\end{defn}

Now, we describe suitable coordinates in an open neighbourhood around a point in the exceptional divisor of $X'$.




\begin{obs} \label{localcoord}
We can construct local coordinates $(\eps,x,v)$ for $X'$ around a point $x'\in E'$  such that:
\begin{enumerate}[i.]
\item $\varepsilon$ is a local single equation defining the Cartier divisor $E'\subset X'$ in an open neighborhood of $x'$.
\item  $(0,x,v)$ represents local coordinates for $E'$ around $x'$ such that $x$ moves along the base $(B_0)_{\text{red}}$ and $v$ moves along the fiber $E'_x$ over $x$ of the projective bundle $E'\simeq \PP N_{(B_0)_{\text{red}}} X$.


\item  the coordinate $x$ represents the projective class of an element also denoted $x\in\Lambda\times (S_n(1))^3$ and the coordinate $v$ represents the projective class of a vector inside $T_{[x]}X\simeq \Lambda\times (S_n(1))^3/\langle x \rangle$ that lifts to a vector also denoted $v\in \Lambda\times (S_n(1))^3$. Then for every $\eps \neq 0$ in a small open neighbourhood of $0$ we have the equality $\pi(\eps,x,v)=[x+\eps v]$. Notice that the expression $[x+\eps v]\in X$ does not depend on the choice of the elements $x,v\in \Lambda\times (S_n(1))^3$.

\end{enumerate}

\end{obs}

With this in mind, we compute the set-theoretical base locus of $\rho'$.

\begin{prop}\label{blprimeraexplosion}
The scheme $B'_0$ is supported on the subvariety $\PP N_{(B_0)_{\text{red}}} B_0 \subset E'$. 
In other words, the reduced scheme $(B'_0)_{\text{red}}$ is equal to $\PP N_{(B_0)_{\text{red}}}B_0$.

\begin{proof}
Let us consider $x' \in E'$ and local coordinates $(\varepsilon,x,v)$ around $x'$ as in Remark \ref{localcoord}. 
Using Lemma \ref{lema1}, for every $\varepsilon\neq 0$ in a small open neighborhood of $x'$ we can compute
\[\rho(\pi(\varepsilon,x,v)) = \rho([x+\varepsilon v]) = \varepsilon H_1(x,v) + \varepsilon^2 H_2(x,v) +\mathcal{O}(\varepsilon^3).\]
Since $H_1(x,v)$ does not vanish generically on $E'$, the generic multiplicity of $\rho([x+\varepsilon v])$ along $E'=(\varepsilon =0)$ is one. As a consequence, we define $\rho'$ locally using the coordinates of Remark \ref{localcoord}:
\begin{equation}\label{rho1}
\rho'(\varepsilon,x,v) := \tfrac{1}{\varepsilon} \rho([x+\varepsilon v]). 
\end{equation}  	
We deduce that $\rho'(0,x,v) = H_1(x,v)$ provided the expression is not equal to zero. Therefore $(B'_0)_{\text{red}}$ consists of the points with coordinates $(0,x,v)$ such that $H_1(x,v)=0$, i.e. such that $v$ represents an element of $\PP N_{(B_0)_{\text{red}},x} B_0$ as claimed.
\end{proof}
\end{prop}

We proceed with the second step of the resolution of $\rho$.

\begin{defn}
Let $X''$ be the blow-up of $X'$ along $(B_0')_{\text{red}}$, let $\pi':X''\rightarrow X'$ be its corresponding projection map and let $E''$ be its exceptional divisor. Also, let $\rho'':X''\dashrightarrow \PP H^0(\PP^n,\Omega_{\PP^n}^1(3))$ be the extension of the map $\rho'$, let $B_{\rho''}$ be its indeterminacy scheme and let $\dbtilde{B_i}$ be the double strict transforms of $B_i$ for $i=1,2,3$.   
\end{defn} 

As in the first step of the resolution of $\rho$, our goal is to calculate the support of the base locus $B_{\rho''}$. For this, we will prove a couple of technical results. But first, we describe suitable coordinates in an open neighbourhood around a point in the exceptional divisor of $X''$.

\begin{obs} \label{localcoord2}

In the same line of Remark \ref{localcoord}, we construct local coordinates $(\eps',x',v')$ in $X''$ around a point $x''\in E''$. Firstly, $\eps'$ is a local equation defining $E''\subseteq X''$. Also, $(0,x',v')$ are local coordinates for $E''$ such that $x'$ moves along the center of the blow-up $(B_0')_{\text{red}}$ and $v'$ moves along the fiber $E''_{x'}$ over $x'$ of the projective bundle $E''\simeq \PP N_{(B_0')_{\text{red}}} X'$. This means that $x'$ represents a point in $X'$ of coordinates $(0,x,v)$ such that $H_1(x,v)=0$ and $v'$ represents the projective class of a normal vector inside $N_{(B'_0)_{\text{red}},x'} X'\simeq T_{x'}X'/T_{x'}(B'_0)_{\text{red}}$. 

Because of Proposition \ref{blprimeraexplosion}, there are $n$ linearly independent directions of the form $[(0,0,\widetilde{v})]$ inside $N_{(B'_0)_{\text{red}},x'} X'$. Adding a direction of the form $[(\widetilde{\eps},0,0)]$, we get $n+1$ linearly independent directions in  $N_{(B'_0)_{\text{red}},x'} X'$. We also know that $\rank(N_{(B_0')_{\text{red}}} X')=n+1$, as $\dim(X'')=3n+1$ and $\dim((B_0')_{\text{red}})=2n$. Therefore, the coordinate $v'$ moves along an element of the form $[(\widetilde{\eps},0,\widetilde{v})]\in X'$ satisfying $H_1(x,\widetilde{v})\neq 0$. Moreover, for every $\eps'\neq 0$ in a small open neighbourhood of $0$ we have the equality $\pi'(\eps',x',v')=(0,x,v)+\eps' (\widetilde{\eps},0,\widetilde{v})$.


\end{obs}

The following Proposition gives an explicit formula for $\rho''$ over the points of the exceptional divisor $E''$ of $X''$.

\begin{prop}
Using the local coordinates of Remark \ref{localcoord2}, we have the equality
	\[\rho''(\varepsilon',x',v') =  H_1(x,\widetilde{v}) + \widetilde{\varepsilon} H_2(x,v) + \mathcal{O}(\varepsilon'),\]	
	provided the expression is not equal to zero.
\begin{proof}
First, following Proposition \ref{blprimeraexplosion}, we compute
	\[\rho'\circ \pi' (\varepsilon',x',v') = \rho'( (0,x,v) + \varepsilon' (\widetilde{\varepsilon},0,\widetilde{v}) ) = \tfrac{1}{\varepsilon' \widetilde{\varepsilon}} \rho([x + \varepsilon'\widetilde{\varepsilon}(v+\varepsilon'\widetilde{v})]). \]
Using Lemma \ref{lema1}, we get
	\begin{align*} \rho'\circ \pi' (\varepsilon',x',v') &= \tfrac{1}{\varepsilon' \widetilde{\varepsilon}} ((\varepsilon')^2 \widetilde{\varepsilon} H_1(x,\widetilde{v}) + (\varepsilon'\widetilde{\varepsilon})^2H_2(x,v) + \mathcal{O}((\varepsilon'\widetilde{\varepsilon})^3))= \\ 
	&= \varepsilon' (H_1(x,\widetilde{v}) + \widetilde{\varepsilon} H_2(x,v)) + \mathcal{O}((\varepsilon')^2).
	\end{align*}
As the expression $H_1(x,\widetilde{v}) + \widetilde{\varepsilon} H_2(x,v)$ does not vanish identically on $E''$, the generic multiplicity of $\rho'\circ \pi'$ along $E''=(\varepsilon' =0)$ is equal to 1. Consequently we define 
\begin{align} \label{rho2}
 \rho'' (\varepsilon',x',v') = \frac{1}{\varepsilon'} \left ( \rho'\circ \pi' (\varepsilon',x',v') \right ) 
 \end{align} 
and we get our claim.
\end{proof}

\end{prop}

In particular, $\rho''$ is well defined over a closed point in the exceptional divisor of local coordinates $(0,(0,x,v),[(\widetilde{\eps},0,\widetilde{v})])$ if and only if $H_1(x,\widetilde{v}) + \widetilde{\varepsilon} H_2(x,v)$ does not vanish. Let us show that it does not vanish in the fiber of a generic point of $(B_0')_{\text{red}}$. For this, we will use the following Lemma.

\begin{lem} \label{lemablsegundaexplosion}
Let $x=[(\cl{\lambda},(F_0,F_0,F_0))]\in (B_0)_{\text{red}}$ be a closed point and let $v=[(\cl{\lambda'},\cl{F'})], \widetilde{v}=[(\cl{\widetilde{\lambda'}},\cl{\widetilde{F'}})]\in T_x X$ with $\cl{F'}, \cl{\widetilde{F'}}\in (S_n(1))^3$. Assume that $H_1(x,v)=H_1(x,\widetilde{v})+\alpha H_2(x,v)=0$ for some $\alpha\in\CC\setminus 0$. Then, either $F'_1\equiv F'_2\equiv F'_3 \text{ mod}(F_0)$ or there exists $i\neq j$ such that $F'_i\equiv F'_j \text{ mod}(F_0)$ and $\lambda_i+\lambda_j=0$.
\begin{proof}
From the equation $H_1(x,v)=d\left(\frac{\sum_{i=1}^3\lambda_i F'_i}{F_0}\right)=0$, we deduce that $\sum_{i=1}^3\lambda_i F'_i\equiv 0 \text{ mod}(F_0)$. Additionally, the equality $H_1(x,\widetilde{v})+\alpha H_2(x,v)=0$ and the fact that $H_1(x,v)=0$ yields
\begin{align*}
d\left(\frac{\sum_{i=1}^3 F_0(\lambda_i\widetilde{F_i'} +\alpha\lambda_i'F'_i)-\alpha\frac{\lambda_i}{2} {F'_i}^2}{F_0^2} \right)=0.
\end{align*}
Thus the expression inside the differential is constant. By multiplying it by $F_0^2$ and reducing it modulo $F_0$, we get $\sum_
{i=1}^3\lambda_i {F'_i}^2\equiv 0 \text{ mod}(F_0)$. Adding the descent condition $\sum_{i=1}^3\lambda_i=0$, we get the following system
\begin{align*}
\begin{pmatrix}
1 & 1 & 1 \\
F'_1 & F'_2 & F'_3 \\
{F'_1}^2 & {F'_2}^2 & {F'_3}^2
\end{pmatrix} \cdot 
\begin{pmatrix}
\lambda_1\\
\lambda_2 \\
\lambda_3
\end{pmatrix}
\equiv 0 \text{ mod}(F_0).
\end{align*}
The result follows directly from a standard determinant and kernel computation of Vandermonde matrices.
\end{proof} 
\end{lem}

We use these last two results to compute the support of $B_{\rho''}$.

\begin{prop} \label{blsegundaexplosion}
The scheme $B_{\rho''}$ is supported in the disjoint union $\bigsqcup_{i=1}^3 \dbtilde{B_i}$.

\begin{proof} 
Observe that $\rho'$ is well-defined outside $(B_0')_{\text{red}}\cup\bigcup_{i=1}^3 B_i'$ and the generic point of the $\dbtilde{B_i}$ belongs to the base locus $B_{\rho''}$. Hence, it is enough to show that the only closed points of $B_{\rho''}$ intersecting the exceptional divisor $E''$ belong to one of the $\dbtilde{B_i}$. Let us take a closed point $x''\in E''$ such that $\rho''(x'')=0$. The local coordinates $(0,(0,x,v),[(\widetilde{\eps},0,\widetilde{v})])$ for $x''$ of Remark \ref{localcoord2}, satisfy the equalities
\begin{align} \label{eqblsegundaexplosion}
\rho(x)=0 & & H_1(x,v)=0 & & H_1(x,\widetilde{v})+\widetilde{\eps}H_2(x,v)=0.
\end{align} 

From these conditions, we deduce that $\widetilde{\eps}\neq 0$. Otherwise, we would get $H_1(x,\widetilde{v})=0$ and therefore $(\widetilde{\eps},0,\widetilde{v})$ would represent the zero vector inside $N_{(B_0')_{\text{red}},x'} X'$. This is a contradiction as $v'$ is the projective class of $(\widetilde{\eps},0,\widetilde{v})$. Notice that $v$ cannot be the class of a vector tangent to $(B_0)_{\text{red}}$ as its projective class inside $N_{(B_0)_{\text{red}},x} X$ should be well defined. Thus, applying Lemma \ref{lemablsegundaexplosion}, using the same notation for $x,v,\widetilde{v}$ and assuming without losing generality that $i=1, j=2$ we get that $F'_1=F'_2\neq F'_3$, $\lambda_1+\lambda_2=\lambda_3=0$. Consequently, we write $x=[((1,-1,0),(F_0,F_0,F_0))]$. Also, we can assume that $\cl{\lambda'}=\vec{0}$ and $F_2'=0$ as $v=[(\cl{\lambda'},\cl{F'})]$ represents a vector of $N_{(B_0)_{\text{red}},x} X$. 
We therefore get that $H_2(x,v)=0$. Also, from equation \eqref{eqblsegundaexplosion} we deduce that $H_1(x,\widetilde{v})=0$, i.e. $v'$ is represented by $[(\widetilde{\eps},0,\widetilde{v})]=[(1,0,0)]$ inside $\PP N_{(B_0')_{\text{red}},x'} X'$. 

Finally, let us prove that $x''=(0,(0,x,v),[(1,0,0)])$ is a closed point of $\dbtilde{B_3}$. For this, we construct locally a curve $x''(t)$ such that $\lim_{t\to 0} x''(t)=x''$. Indeed, the curve given in local coordinates of Remark \ref{localcoord2} by $x''(t)=(t,(0,x,v),[(1,0,0)])$ pushes via the two blow-up projections to the curve 
\begin{align*}
\pi\circ\pi'(x''(t))=\pi((0,x,v)+t(1,0,0))=[x+tv]=[((1,-1,0),(F_0,F_0,F_0+tF'_3))].
\end{align*}
But $\pi\circ\pi'(x''(t))$ is contained in $B_3$ for every $t\neq 0$. The proposition follows.
\end{proof}

\end{prop}

We arrive at the third and final step of the resolution of $\rho$. 

\begin{defn}
Let $X'''$ be the blow-up of $X''$ along $\bigsqcup_{i=1}^3 \dbtilde{B_i}$, let $\pi'':X'''\rightarrow X''$ be its projection map and let $\rho''':X'''\dashrightarrow \PP H^0(\PP^n,\Omega_{\PP^n}^1(3))$ be the extension of $\rho''$. Also, let $B_{\rho'''}$ be its base locus scheme, let $E'''_i=\pi''^{-1}(\dbtilde{B_i})$ be the fiber of $\dbtilde{B_i}$ via $\pi''$ and let $E'''=\bigsqcup_{i=1}^3 E'''_i$ be the exceptional divisor. 
\end{defn}

Before showing that the map $\rho'''$ is regular, let us describe suitable coordinates in an open neighbourhood around a point in the exceptional divisor of $X''$.

\begin{obs} \label{localcoord3}

Using Remarks \ref{localcoord} and \ref{localcoord2}, we construct local coordinates $(\eps'',x'',v'')$ in $X''$ around a point $x'''\in E'''$ in a similar fashion. Here $\eps''$ is a local equation for $E'''\subseteq X'''$. As $B_{\rho''}$ is smooth along $\dbtilde{B_i}\setminus \left(\dbtilde{B_i}\cap E''\right)$, we will focus on the points $x'''\in \pi''^{-1}\left(E'''\cap \dbtilde{B_i}\right)$ and assume that $i=3$. 


The coordinate $x''$ represents a point $[(0,x',v')]\in X''$. Here $x'$ is of the form $[(0,x,v)]\in X'$ such that $H_1(x,v)=0$, $x=[((1,-1,0),(F_0,F_0,F_0))]\in (B_0)_{\text{red}}$, $v=(\vec{0},(0,0,F'_3))\in \Lambda \times \left(S_n(1)\right)^3$ and $v'=[(1,0,0)]$.

 For the coordinate $v''$, we find a basis of $n+1$ vectors of $N_{B''_3,x''}X''$. We select $n$ of them with coordinates $(0,0,\widetilde{v}')\in T_{x''}X''$, where $\widetilde{v}'=(0,0,\dobletilde{v})\in T_{x'}X'$. Because of Remark \ref{secB0}, $\dobletilde{v}$ represents a direction normal to $B_3$ such that $H_1(x,\dobletilde{v})=0$ and thus it lifts to a vector of the form $(\vec{0},(\minidobletilde{F_1'},0,0))\in \Lambda \times \left(S_n(1)\right)^3$. We add a direction normal to $B_3''$ of the form $(0,\widetilde{x}',0)\in T_{x''}X''$ where $\widetilde{x}'$ has coordinates $(0,\widetilde{x},0)\in T_{x'}X'$ and $\widetilde{x}=((0,1,-1),\vec{0})\in \Lambda \times \left(S_n(1)\right)^3$. Therefore, a general vector has the form $v''=(\vec{0},\widetilde{x}',\widetilde{v}')$ where $\widetilde{x}'=(0,\alpha\widetilde{x},0)\in X'$ for some $\alpha\in\CC$ and $\widetilde{v}'=(0,0,\dobletilde{v})\in T_{x'}X'$ as before.


\end{obs}

To finish the section, let us show that the induced map $\rho''':X'''\rightarrow \PP H^0(\PP^n,\Omega_{\PP^n}^1(3))$ with base locus scheme $B_{\rho'''}$ is indeed regular.

\begin{teo}
The map $\rho''':X'''\rightarrow \PP H^0(\PP^n,\Omega_{\PP^n}^1(3))$ depicted in Diagram \ref{diag} is a resolution of the rational parametrization map $\rho:X\dashrightarrow \PP H^0(\PP^n,\Omega_{\PP^n}^1(3))$ of the logarithmic component $\LL(1,1,1)$ of type $(1,1,1)$.


\begin{proof}
We only need to prove that the induced map $\rho'''$ is regular. Without loss of generality, it is enough to prove that $\rho'''$ is well defined in the fiber of a point $x''\in E''\cap\dbtilde{B_3}$. Take local coordinates $(\eps'',x'',v'')$ around $x'''$ as in Remark \ref{localcoord3}. Now, we use Lemma \ref{lema1} and the definition of $\rho''$ (equation \eqref{rho2}) to get a formula for $\rho'''$. We get

\begin{align} \label{rho3}
\begin{aligned}
\rho''\circ\pi''(\eps'',x'',v'')=\rho''(0,x'+\eps''\widetilde{x'},v'+\eps''\widetilde{v'})=H_1(x+\eps''\alpha\widetilde{x},\eps''\dobletilde{v})+H_2(x+\eps''\alpha\widetilde{x},v)= \\
\eps''(H_1(x,\dobletilde{v})+H_2(\widetilde{x}_{\alpha,F_0},v))+\mathcal{O}(\eps''^2)
\end{aligned}
\end{align}
where $\widetilde{x}_{\alpha,F_0}=((0,\alpha,-\alpha),(F_0,F_0,F_0))$. Cancelling out $\eps''$ in \eqref{rho3} we get the formula for $\rho'''$. In particular, inside $E'''$ it is given by the formula 
\begin{align} \label{absfinal}
\rho'''(0,x'',v'')=H_1(x,\dobletilde{v})+H_2(\widetilde{x}_{\alpha,F_0},v). 
\end{align}
Now, suppose there exists a point of coordinates $(0,x'',v'')$ that belongs to the base locus $B_{\rho'''}$ of $\rho'''$. Then, the expression for $\rho'''$ in equation \eqref{absfinal} vanishes. Dividing by $F_0^3$ we get
\begin{align*}
d\left(\frac{\dobletilde{F'_1}}{F_0}\right)=0.
\end{align*} 
This yields $\dobletilde{F'_1}\equiv 0 \text{ mod} (F_0)$ and therefore $\dobletilde{v}$ is zero. This is a contradiction because $\dobletilde{v}$ cannot vanish, as $v''$ is the projective class of a normal vector.
\end{proof}

\end{teo}

The following commutative diagram summarizes this section: 


\begin{center}
\begin{tikzcd} 
 & E'''=\bigsqcup_{i=1}^3 E_i''' \arrow{d} \arrow[hookrightarrow]{r} & X'''\arrow{dddr}{\rho'''} \arrow{d}{\pi''} & \\
 E''\arrow{d} \arrow[hookrightarrow,bend right=20,swap]{rr} & \bigsqcup_{i=1}^3 \dbtilde{B_i} \arrow[hookrightarrow]{r} & X'' \arrow[dashed]{ddr}{\rho''} \arrow{d}{\pi'} & \\
(B_0')_{\text{red}} \arrow[hookrightarrow]{r} & E'\arrow{d} \arrow[hookrightarrow]{r}  & X' \arrow{d}{\pi} \arrow[dashed]{dr}{\rho'} & \\
& (B_0)_{\text{red}} \arrow[hookrightarrow]{r} & X \arrow[dashed]{r}{\rho} & \PP H^0(\PP^n,\Omega_{\PP^n}^1(3))
\end{tikzcd}
\captionof{Diagram}{The sequence of blow-ups resolving the rational parametrization map $\rho$.}
\label{diag}
\end{center}

\section{Degree Calculations}

 In this section, we use the resolution of the rational parametrization map $\rho$ of Section 3 compute the degree of the logarithmic component of type $(1,1,1)$. For this, we use standard intersection theory techniques. 
 
 Let $h$ be the hyperplane class of $\PP H^0(\PP^n ,\Omega_{\PP^n}^1(3))$. We know that the degree of $\LL(1,1,1)$ is defined as the degree of the Chow 0-cycle $(h|_{\LL(1,1,1)})^{3n+1}$ inside $\LL(1,1,1)$. Recall that the map $\rho$ is generically injective up to order, see Proposition \ref{geninyupto}. Thus, the degree of the map $\rho$ is the same as the order of the group $\mathbb{S}(1,1,1)\simeq \mathbb{S}_3$ defined in Remark \ref{defS(d)}. We deduce that the map $\rho'''$ is generically 6 to 1 and 
 \begin{align*}
 \int_{X'''}\rho'''^*\left(h\right)^{3n+1}=6\int_{X'''}(h|_{\LL(1,1,1)})^{3n+1}. 
 \end{align*}
 Therefore, to compute the degree of $\LL(1,1,1)$, we pullback the hyperplane class $h$ through $\rho'''$. Recall that the rational map $\rho$ was defined by sections of the line bundle $\mathcal{O}_{X}(1,1,1,1)$. Using equations \eqref{rho1}, \eqref{rho2} and \eqref{rho3} and the local calculations in the previous section, we know that $\rho'''$ is defined by sections of the line bundle 
$$\mathcal{O}_{X}(1,1,1,1) \otimes \mathcal{O}_{X'}(-E')\otimes \mathcal{O}_{X''}(-E'') \otimes \oplus_{i=1}^3\mathcal{O}_{X'''}(-E_i''').$$
 We conclude that the pullback of the hyperplane class $h$ of the target space of $\rho'''$ is  
\begin{align*} 
\rho'''^{*}(h)=\sum_{i=1}^4h_i-e_1-e_2-\sum_{i=1}^3e_{3,i}.
\end{align*}
Here $h_1,h_2,h_3,h_4$ are generators of $X$ such that $h_1$ is the (pullback of the) hyperplane class of $\PP\Lambda\simeq\PP^1$ and $h_2,h_3,h_4$ are the (pullback of the) hyperplane class of the three different $\PP S_n(1)\simeq\PP^n$. Also, $e_1=[E']$, $e_2=[E'']$ and $e_{3,i}=[E'''_{i}]$.
Consequently, the degree of $\LL(1,1,1)$ is computed as 
\begin{align*}
\frac{1}{6}\int_{X'''}\left(\sum_{i=1}^4h_i-e_1-e_2-\sum_{i=1}^3e_{3,i}\right)^{3n+1}.
\end{align*} 

For computing this number, we use standard intersection theory techniques. The main reference for this theory is \cite{F}. Using the fact that the three blow-up projections are birational, we will push forward the zero cycle $\rho'''^*(h)^{3n+1}$ to $X$ in three steps to get a zero cycle inside $X$ with the same degree. The computation inside $X$ is easier as it is a product of projective spaces. 

For pushing forward this cycle, we need to understand the pushforwards of powers of the cycles $e_1,e_2,e_{3,i}$ through their corresponding blow-up projections. These are exactly the \textit{Segre classes} of the normal bundles of the respective centers of blow-up in the ambient space, see \cite[Corollary 4.2.2]{F}. Eventhough Segre classes can be defined generally for every cone inside a scheme, in the case of vector bundles they can be defined as the graded parts of a formal inverse of the corresponding total Chern class inside the Chow ring, see \cite[Proposition 4.1]{F}. The sum of all the Segre classes of a vector bundle is the \textit{total Segre Class} and it is additive in exact sequences. 
   
As a first step, we collect the coefficients of $e_{3,i}$ and use our knowledge of the normal bundle $N_{\dbtilde{B_i}} X''$ to push it forward to $X''$. Indeed, we know that 
\begin{align*}
\pi''_*(e_{3,i}^j)=(-1)^{j-1}s_{j'}(N_{\dbtilde{B_i}} X'')\cap [\dbtilde{B_i}]
\end{align*}
where $j'=j-\text{codim}(\dbtilde{B_i})=j-n-1$. 

For computing the class $[\dbtilde{B_i}]$ we use twice the Blow-Up Formula, see \cite[Theorem 6.7]{F}:
\begin{align*}
[\dbtilde{B_i}]=[B_i]-{j_1}_*\left\{c(\EE_1)\cap \pi^*(s(N_{B_i\cap (B_0)_{\text{red}}}{B_i})\right\}_{n+1} \\
- {j_2}_*\left\{c(\EE_2)\cap \pi'^*(s(N_{B_i'\cap (B_0')_{\text{red}}}B_i'))\right\}_{n+1}.
\end{align*}
Here $j_1$ and $j_2$ are the inclusions of $E'$ and $E''$ in $X'$ and $X''$ and the bundles $\EE_1\simeq \pi^*(N_{(B_0)_{\text{red}}} X)/N_{E'}X'$ and $\EE_2\simeq\pi'^*(N_{(B_0')_{\text{red}}} X')/N_{E''}X''$ are the \textit{excess bundles}. The Segre class $s(N_{B_i\cap (B_0)_{\text{red}}}{B_i})$ is deduced from the exact sequence
\begin{align*}
0\rightarrow T(B_i\cap (B_0)_{\text{red}})\rightarrow T B_i|_{B_i\cap (B_0)_{\text{red}}}\rightarrow N_{B_i\cap (B_0)_{\text{red}}}B_i \rightarrow 0. 
\end{align*} 
On the other hand, we deduce the equality $B_i'\cap (B_0')_{\text{red}}=B_i'\cap E'$ using the description of the tangent space $T_{B_0}$ in equation \eqref{eq3}. Therefore, we get $N_{B_i'\cap (B_0')_{\text{red}}}B_i'\simeq \mathcal{O}_{E'_{B_i'}}(-1)$ where $E'_{B_i'}$ is the exceptional divisor of $B_i'$. For the computation of $c(E_1)$  we use the isomorphism  $N_{E'}X'\simeq \mathcal{O}_{E'}(-1)$ and the exact sequence 

\begin{align} \label{NB0redX}
0\rightarrow T(B_0)_{\text{red}}\rightarrow T X|_{(B_0)_{\text{red}}}\rightarrow N_{(B_0)_{\text{red}}}X\rightarrow 0.
\end{align}
To calculate $c(E_2)$ we use the two isomorphisms $N_{(B_0')_{\text{red}}}E'\simeq \mathcal{O}_{E'}(1)\otimes \pi'^*(N_{(B_0)_{\text{red}}}X/N_{(B_0)_{\text{red}}}B_0)$ and $N_{E''}X''\simeq \mathcal{O}_{E''}(-1)$, the exact sequence

\begin{align} \label{NB0'redX'}
0\rightarrow N_{(B'_0)_{\text{red}}} E'\rightarrow N_{(B'_0)_{\text{red}}} X'\rightarrow N_{E'} X'|_{(B'_0)_{\text{red}}}\rightarrow 0
\end{align} 
 and the formula for the Chern class of a tensor product with a line bundle , see \cite[Proposition 9.13]{EH} and \cite[Proposition 5.17]{EH}. The Chern classes of the bundle $N_{(B_0)_{\text{red}}}X/N_{(B_0)_{\text{red}}}B_0$ can be deduced from the exact sequences of Remark \ref{secB0}, the exact sequences
\begin{align*}
0 \rightarrow T (B_0)_{\text{red}} \rightarrow T X|_{(B_0)_{\text{red}}} \rightarrow N_{(B_0)_{\text{red}}} X \rightarrow 0 \\
0 \rightarrow T (B_0)_{\text{red}} \rightarrow T B_0|_{(B_0)_{\text{red}}} \rightarrow N_{(B_0)_{\text{red}}} B_0 \rightarrow 0 
\end{align*} 
and our knowledge of the tangent bundles of $(B_0)_{\text{red}}$ and $X$. 

For the calculation of the Segre classes of the normal bundles $N_{\dbtilde{B_i}}X''$, consider the auxiliary algebraic subsets of $X$ of dimension $2n+1$ containing both $(B_0)_{\text{red}}$ and the corresponding $B_i$
$$\begin{cases}
Z_1=\left\{[(\cl{\lambda},(A,B,B))]: A,B\in S_n(1)\right\}\subseteq X  \\
Z_2=\left\{[(\cl{\lambda},(B,A,B))]: A,B\in S_n(1)\right\}\subseteq X  \\
Z_3=\left\{[(\cl{\lambda},(B,B,A))]: A,B\in S_n(1)\right\}\subseteq X.
\end{cases}$$
The key fact is that $(B_0)_{\text{red}}$ and $B_i$ intersect properly inside $Z_i$. This follows from the dimensions and set theoretical descriptions of the sets involved. Thus, we get $s(N_{\dbtilde{B_i}}X'')=s(N_{B_i'}Z_i')s(\pi'^*(N_{Z_i'} X')\otimes \mathcal{O}_{E''}(1))$, where $Z_i'$ is the strict transform of $Z_i$, see \cite[Section 4.3]{Aluffi}. Notice that $B_i$ and $(B_0)_{\text{red}}$ intersect properly inside $Z_i$ and therefore $s(N_{B_i'}Z_i')=s(\pi^*N_{B_i}Z_i)$. From the exact sequence 
$$0\rightarrow T B_i\rightarrow T Z_i|_{B_i}\rightarrow N_{B_i}Z_i\rightarrow 0,$$
we get that $s(\pi^*N_{B_i}Z_i)=\frac{1}{1+h_1}|_{B_i}=1$. For computing $s(\pi'^*(N_{Z_i'} X')\otimes \mathcal{O}_{E''}(1))$ we use the isomorphism $N_{Z_i'}X'\simeq N_{Z_i}X \otimes \mathcal{O}_{E'}(1)$ (see \cite[Appendix B.6.10]{F}), the formula for the Segre class of a tensor product with a line bundle (\cite[Example 3.1.1]{F}) and the exact sequence 
$$0\rightarrow T Z_i\rightarrow T X|_{Z_i}\rightarrow N_{Z_i}X\rightarrow 0.$$

As a second step, we collect the terms of coefficient $e_2^j$. We know that 
$$\pi'_*(e_2^j)=(-1)^{j-1}s_{j'}(N_{(B_0')_{\text{red}}}X')\cap [(B_0')_{\text{red}}]$$
where $j'=j-\text{codim}((B_0')_{\text{red}})=j-n-1$. The computation of the class $[(B_0')_{\text{red}}]$ is deduced from the formula for the class of a projective bundle in terms of the hyperplane class of the exceptional divisor and the already mentioned Chern classes of the bundle $\pi'^*(N_{(B_0)_{\text{red}}} X/N_{(B_0)_{\text{red}}} B_0)$, see \cite[Proposition 9.13]{EH}. For the Segre classes $s_{j'}(N_{(B_0')_{\text{red}}}X')$, we use the isomorphisms $N_{(B_0')_{\text{red}}}E'\simeq \mathcal{O}_{E'}(1)\otimes \pi'^*(N_{(B_0)_{\text{red}}}X/N_{(B_0)_{\text{red}}}B_0)$ and $N_{E'} X'\simeq \mathcal{O}_{E'}(-1)$, the exact sequence \eqref{NB0'redX'} and the formula for the Segre class of a tensor product with a line bundle, see \cite[Example 3.1.1]{F}. 

Finally, we collect the terms of coefficient $e_1^j$. We know that 
$$\pi_*(e_1^j)=(-1)^{j-1}s_{j'}(N_{(B_0)_{\text{red}}}X)\cap [(B_0)_{\text{red}}],$$
where $j'=j-\text{codim}((B_0')_{\text{red}})=j-2n$. As $(B_0)_{\text{red}}$ is the image of a diagonal embedding, its class inside the Chow ring of $X$ is $[(B_0)_{\text{red}}]=\sum_{0\leq i,j\leq n}h_2^ih_3^jh_4^{2n-i-j}$. For the Segre classes $s_{j'}(N_{(B_0)_{\text{red}}} X)$ we use the exact sequence \eqref{NB0redX}.

To implement the actual calculations, it is best to use a computer program. A script using Macaulay2 can be found in \cite{Programa}. 
The following table shows the degree of the logarithmic components $\LL(1,1,1)$ of type $(1,1,1)$ in $\PP^n$ for small $n$.

\begin{center}
\begin{tabular}{ |c|c| } 
\hline
n  & Degree \\
\hline
3  & 80 \\ 
\hline
4  & 4035 \\
\hline
5  & 165984 \\
\hline
6  & 6091960 \\
\hline
7  & 208063680 \\
\hline
8  & 6766823415 \\
\hline 
\end{tabular}
\captionof{table}{Degree of the logarithmic components $\LL(1,1,1)$ of type $(1,1,1)$ in $\PP^n$ for small $n$.}
\end{center}


In particular, we have the following

\begin{teo} \label{grado111}
The degree of the logarithmic component $\LL(1,1,1)$ of type $(1,1,1)$ in $\PP^3$ is $80$.
\end{teo}

The degree of the component $\LL(2,1)$ inside $\mathbb{F}(3,1)$ is 55. It was computed in \cite[Section 5.3.1]{CPV}. Therefore, this completes the calculation of the degree of the two irreducible components of the space of Jacobi equations in $\PP^3$. 

\begin{cor}
The moduli space $\mathbb{F}(3,1)$ of Jacobi equations in $\PP^3$ has two irreducible components $\LL(2,1)$ and $\LL(1,1,1)$ of degrees $55$ and $80$ respectively. 

\end{cor} 




\vspace{0.5cm}

$^1$ Departamento de Matemática - IMAS, FCEyN, Universidad de Buenos Aires, Buenos Aires, Argentina. The author was fully supported by CONICET.

\textit{E-mail address: }\href{mailto:mchehebar@dm.uba.ar}{mchehebar@dm.uba.ar}


\bibliographystyle{plain}
\bibliography{template}

\end{document}